\documentclass[oneside,brazil,english]{amsart}
\usepackage[T1]{fontenc}
\usepackage[latin9]{inputenc}
\usepackage{geometry}
\usepackage{amstext}
\usepackage{amsthm}
\usepackage{amssymb}

\makeatletter
\numberwithin{equation}{section}
\numberwithin{figure}{section}
  \theoremstyle{plain}
  \newtheorem*{thm*}{\protect\theoremname}

\@ifundefined{date}{}{\date{}}
\usepackage{amsmath}
\usepackage{microtype}

\makeatother

\usepackage{babel}
  \addto\captionsbrazil{\renewcommand{\theoremname}{Teorema}}
  \addto\captionsenglish{\renewcommand{\theoremname}{Theorem}}
  \providecommand{\theoremname}{Theorem}

\begin{document}

\title{On the number of roots of self-inversive polynomials on the complex
unit circle}

\author{R. S. Vieira}

\email{rsvieira@df.ufscar.br}

\maketitle
\begin{center}Universidade Federal de São Carlos, Departamento de
Física, CEP 13560-905, São Carlos, Brazil.\end{center}
\begin{abstract}
We present a sufficient condition for a self-inversive polynomial
to have a fixed number of roots on the complex unit circle. We also
prove that these roots are simple when that condition is satisfied.
This generalizes the condition found by Lakatos and Losonczi for all
the roots of a self-inversive polynomial to lie on the complex unit
circle.

\vspace{0.5cm}


\vspace{1.5cm}

\setcounter{section}{1}
\end{abstract}

A polynomial 
\begin{equation}
p(z)=a_{0}+a_{1}z+...+a_{n-1}z^{n-1}+a_{n}z^{n},\label{Pz}
\end{equation}
with coefficients in $\mathbb{C}$ and with $a_{n}\neq0$ is called
a $n$th degree \emph{self-inversive polynomial} if it satisfies the
property 
\begin{equation}
p(z)=\omega z^{n}\overline{p}\left(\frac{1}{z}\right),\qquad\text{with}\qquad\left|\omega\right|=1,
\end{equation}
where $\overline{p}(z)$ is the complex-conjugate of $p(z)$. In the
special case $\omega=1$, $p(z)$ is called a \emph{self-reciprocal
polynomial}. From this definition follows that the coefficients of
a self-inversive polynomial satisfy the relations 
\begin{equation}
a_{n-k}=\omega\overline{a}_{k},\qquad0\leq k\leq n,\label{Ak0}
\end{equation}
so that 
\begin{equation}
\left|a_{n-k}\right|=\left|a_{k}\right|,\qquad0\leq k\leq n.\label{Ak}
\end{equation}

In \cite{Lakatos} Lakatos and Losonczi presented sufficient conditions
for all the roots of a self-inversive polynomial $p(z)$ to lie on
the complex unit circle $U=\left\{ z\in C:\left|z\right|=1\right\} $.
Here we shall extend this theorem by presenting sufficient conditions
for a self-inversive polynomial to have a fixed number of roots on
$U$. Namely, we shall prove the following
\begin{thm*}
Let $p(z)=a_{n}z^{n}+a_{n-1}z^{n-1}+...+a_{1}z+a_{0}$ be a $n$th
degree self-inversive polynomial. If the inequality
\begin{equation}
\left|a_{n-l}\right|>\frac{1}{2}\left(\frac{n}{n-2l}\right)\sum_{\begin{smallmatrix}k=0\\
k\neq\{l,n-l\}
\end{smallmatrix}}^{n}\left|a_{k}\right|,\qquad l<n/2,\label{Anl}
\end{equation}
holds then $p(z)$ has exactly $n-2l$ roots on the complex unit circle
$U$ and these roots are simple. Moreover, if $n$ is even and $l=n/2$,
then $p(z)$ has no roots on $U$ if the inequality 
\begin{equation}
\left|a_{n/2}\right|>\sum_{\begin{smallmatrix}k=0\\
k\neq n/2
\end{smallmatrix}}^{n}\left|a_{k}\right|\label{an2}
\end{equation}
is satisfied.
\end{thm*}
In order to prove this theorem we shall use only two theorems of complex
analysis, namely, the theorems of Cohn \cite{Cohn} and Rouché \cite{Rouche}.
Cohn's theorem states that a self-inversive polynomial $p(z)$ has
as many roots inside $U$ as has the polynomial 
\begin{equation}
q(z)=z^{n-1}\,\overline{p'}\left(\frac{1}{z}\right)\label{Qz}
\end{equation}
in the same region, where $p'(z)$ is the formal derivative of $p(z)$
\cite{Bonsall,Ancochea,Schinzel}. The Rouché theorem states that
if two functions $f(z)$ and $g(z)$ are analytic inside a simple
closed Jordan curve $\gamma$, are continuous on $\gamma$ and $\left|f(z)\right|>\left|g(z)\right|$
in all points of the curve $\gamma$ then $f(z)$ and $h(z)=f(z)+g(z)$
have the same number of roots inside the region delimited by the curve
$\gamma$ \cite{Marden}.
\begin{proof}
Let $p(z)$ be a $n$th degree self-inversive polynomial as in (\ref{Pz}).
From (\ref{Qz}) the $q(z)$ polynomial used in Cohn's theorem is
given by 
\begin{equation}
q(z)=n\bar{a}_{n}+(n-1)\bar{a}_{n-1}z+(n-2)\bar{a}_{n-2}z^{2}+...+2\bar{a}_{2}z^{n-2}+\bar{a}_{1}z^{n-1}.
\end{equation}
Supposing further that $a_{n-l}\neq0$ and $l<n/2$, let us define
the functions 
\begin{equation}
f(z)=(n-l)\bar{a}_{n-l}z^{l},
\end{equation}
 and 
\begin{eqnarray*}
g(z) & = & q(z)-f(z),\\
 & = & n\bar{a}_{n}+...+(n-l+1)\bar{a}_{n-l+1}z^{l-1}+(n-l-1)\bar{a}_{n-l-1}z^{l+1}+...+\bar{a}_{1}z^{n-1}.
\end{eqnarray*}
On the curve $\gamma=\{z\in\mathbb{C}:z=e^{iy},y\in\mathbb{R}\}$,
we have 
\begin{equation}
\left|f(z)\right|_{\gamma}=(n-l)\left|a_{n-l}\right|,
\end{equation}
 and 
\begin{eqnarray}
\left|g(z)\right|_{\gamma} & = & \left|n\bar{a}_{n}+(n-1)\bar{a}_{n-1}e^{iy}+...+(n-l+1)\bar{a}_{n-l+1}e^{iy(l-1)}\right.\nonumber \\
 & + & \left.(n-l-1)\bar{a}_{n-l-1}e^{iy(l+1)}+...+2\bar{a}_{2}e^{iy(n-2)}+\bar{a}_{1}e^{iy(n-1)}\right|,
\end{eqnarray}
but, from the triangle inequality, follows that 
\begin{eqnarray}
\left|g(z)\right|_{\gamma} & \leq & n\left|a_{n}\right|+...+(n-l+1)\left|a_{n-l+1}\right|+(n-l-1)\left|a_{n-l-1}\right|+...+\left|a_{1}\right|.\nonumber \\
\end{eqnarray}
Therefore, the theorem of Rouché can be applied whenever the condition
\begin{equation}
(n-l)\left|a_{n-l}\right|>n\left|a_{n}\right|+...+(n-l+1)\left|a_{n-l+1}\right|+(n-l-1)\left|a_{n-l-1}\right|+...+\left|a_{1}\right|\label{In}
\end{equation}
holds. Moreover, since $p(z)$ is self-inversive, we can use the property
(\ref{Ak}), noticing that 
\begin{equation}
k\left|a_{k}\right|+(n-k)\left|a_{n-k}\right|=n\left|a_{k}\right|=\frac{n}{2}\left(\left|a_{k}\right|+\left|a_{n-k}\right|\right),
\end{equation}
in order to rewrite (\ref{In}) as 
\begin{equation}
(n-l)\left|a_{n-l}\right|>l\left|a_{n-l}\right|+\frac{n}{2}\sum_{\begin{smallmatrix}k=0\\
k\neq\{l,n-l\}
\end{smallmatrix}}^{n}\left|a_{k}\right|.
\end{equation}
Then, the condition for the applicability of the Rouché theorem becomes
\begin{equation}
\left|a_{n-l}\right|>\frac{1}{2}\left(\frac{n}{n-2l}\right)\sum_{\begin{smallmatrix}k=0\\
k\neq\{l,n-l\}
\end{smallmatrix}}^{n}\left|a_{k}\right|.
\end{equation}
When this condition is fulfilled, we get from the Rouché theorem that
$f(z)$ and $q(z)=f(z)+g(z)$ must have the same number of roots inside
the complex unit circle. But $f(z)$ has exactly $l$ roots inside
$U$ and so also has $q(z)$. On the other hand, from the theorem
of Cohn we get that $p(z)$ has $l$ roots inside $U$ as well, but
$p(z)$ is self-inversive and, thus, $p(z)$ also has $l$ roots outside
$U$. Therefore we conclude that $p(z)$ has $n-2l$ roots on the
complex unit circle $U$. 

To prove that $p(z)$ has no multiple roots on $U$ we can proceed
in a very standard way, but now, instead to use the polynomial $q(z)$
we should take the polynomial $p'(z)$, the formal derivative of $p(z)$,
\begin{equation}
p'(z)=a_{1}+2a_{2}z+...+(n-1)a_{n-1}z^{n-2}+na_{n}z^{n-1}.
\end{equation}
Then we define the functions functions $f(z)$ and $g(z)$ used in
Rouché's theorem as 
\begin{equation}
f(z)=(n-l)a_{n-l}z^{n-l-1},\qquad g(z)=p'(z)-f(z).
\end{equation}
From the self-inversive property of $p(z)$ it follows that the same
inequality (\ref{Anl}) holds here, for the validity of the Rouché
theorem in the region delimited by the curve $\gamma$. This means
that $f(z)$ and $p'(z)$ have the same number of roots inside the
complex unit circle, that is, $n-l-1$ roots. But we have seen that
the polynomial $q(z)$ has exactly $l$ roots inside $U$ and, therefore,
$p'(z)$ must have $l$ roots outside $U$ as well (because $q(z)$
is proportional to $\overline{p'}(1/z)$). Since the degree of $p'(z)$
is $n-1$, it follows from the Fundamental Theorem of Algebra that
$p'(z)$ has no roots on $U$. This implies that $p(z)$ has no multiple
roots on $U$ since, if it had, then $p(z)$ would necessarily share
any root with $p'(z)$, leading to a contradiction.

Finally, to prove that no roots of $p(z)$ are on $U$ when the inequality
(\ref{an2}) is satisfied, we just need to define $f(z)=\left(n/2\right)\bar{a}_{n/2}z^{n/2}$
and $g(z)=q(z)-f(z)$ and following the same steps as above. 
\end{proof}
When the condition (\ref{Anl}) is satisfied, the theorem proved above
ensures that a $n$th degree self-inversive polynomial $p(z)=a_{n}z^{n}+...+a_{0}$
has \emph{exactly} $n-2l$ roots on the complex unit circle, $U$.
Nevertheless, a weaker condition than (\ref{Anl}) can also be deduced,
as was indicated by the referee of this paper. Namely, we can show
that a $n$th degree self-inversive polynomial $p(z)$ should have
\emph{at least} $n-2l$ roots on $U$ if the inequality 
\begin{equation}
\left|a_{n-l}\right|>\frac{1}{2}\sum_{\begin{smallmatrix}k=0\\
k\neq\{l,n-l\}
\end{smallmatrix}}^{n}\left|a_{k}\right|\label{Anl2}
\end{equation}
is verified. This result can be shown as follows: first notice that
on the complex unit circle $U$ the polynomial $r(z)=\omega^{-1/2}z^{-n/2}p(z)$
is a real trigonometric polynomial. In fact this can be easily seen
through the relations (\ref{Ak0}) and introducing the change of variables
\begin{equation}
z=e^{it},\qquad\omega=e^{i\sigma},\qquad\mbox{and}\qquad a_{k}=\left|a_{k}\right|e^{i\phi_{k}},\qquad0\leq k\leq n.
\end{equation}
In particular, we have that 
\begin{equation}
\omega^{-1/2}z^{-n/2}\left(a_{l}z^{l}+a_{n-l}z^{n-l}\right)=2\left|a_{n-l}\right|\cos\left[\left(\frac{n}{2}-l\right)t+\phi_{n-l}-\frac{\sigma}{2}\right].
\end{equation}
From this we can see that $\left|a_{l}z^{l}+a_{n-l}z^{n-l}\right|\leq2\left|a_{n-l}\right|$
on $U$, with equality happening only when $t=(2\pi j+\sigma-2\phi_{n-l})/(n-2l)$
for integer $j$. Thus, the validity of the (\ref{Anl2}) means that
the paired term $a_{l}z^{l}+a_{n-l}z^{n-l}$ dominates all the other
terms of $p(z)$. Since $r(z)$ and $\cos\left[\left(n/2-l\right)t+\phi_{n-l}-\sigma/2\right]$
has the same sign on those values of $t$, it follows by the Intermediate
Value Theorem that $r(z)$, and consequently $p(z)$, has \emph{at
least} one root in any interval of $t$ delimited by the successive
values of $j$ on the formula above. This implies that $r(z)$, and
hence $p(z)$, has\emph{ at least} $n-2l$ roots on $U$ when (\ref{Anl2})
is satisfied.

We highlight however that the conditions (\ref{Anl}) and (\ref{Anl2})
are not equivalent, since there exist self-inversive polynomials which
satisfy (\ref{Anl2}) but have more than $n-2l$ roots on $U$. For
instance, the polynomial 
\begin{equation}
h(z)=(1+i)z^{4}-2iz^{3}-2z+(1+i)
\end{equation}
is self-inversive with $\omega=i$, satisfies (\ref{Anl2}) for $l=1$,
but all roots of $h(z)$ are on $U$. 

Now let us comment about some special cases of (\ref{Anl}). First,
setting $l=0$ in the equation (\ref{Anl}), we get the condition
\begin{equation}
\left|a_{n}\right|>\frac{1}{2}\sum_{k=1}^{n-1}\left|a_{k}\right|,\label{Lakatos}
\end{equation}
for all the roots of a $n$th degree self-inversive polynomial to
lie on the complex unit circle $U$. This condition was in fact found
before by Lakatos and Losonczi in \cite{Lakatos}. Lakatos and Losonczi
also showed that the roots of $p(z)$ are simple except when an equality
takes place in the formula (\ref{Lakatos}), in which case the polynomial
may have a double root on $U$ $-$ the authors also gave the exact
conditions for that happen \cite{Lakatos}. It is clear that a self-inversive
polynomial may present multiple roots on $U$ when an equality occurs
in (\ref{Anl}), but the exact conditions to that happen should be
more analyzed.

The case $l=1$ is also interesting, since from the condition 
\begin{equation}
\left|a_{n-1}\right|>\frac{1}{2}\left(\frac{n}{n-2}\right)\sum_{\begin{smallmatrix}k=0\\
k\neq\{1,n-1\}
\end{smallmatrix}}^{n}\left|a_{k}\right|\label{Salem}
\end{equation}
we can test if a given polynomial with integer coefficients is a Salem
Polynomial, as well as we can construct Salem polynomials by giving
to the coefficient $a_{n-1}$ of a self-reciprocal polynomial a larger
enough value. A Salem polynomial is a self-reciprocal polynomial with
integer coefficients whose roots lie all on the complex unit circle,
except for two positive and reciprocal roots $r$ and $1/r$ \cite{Smyth,Boyd2,Boyd,Hironaka}.
Notice however that the condition (\ref{Salem}) is sufficient but
not necessary. Thus, it is not mandatory that a Salem polynomial satisfies
(\ref{Salem}) and, indeed, Salem polynomials with small Salem numbers
\cite{Boyd2,Boyd} known to date do not satisfy this inequality.

This study was born in the research of Vieira and Lima-Santos on the
solutions of the Bethe Ansatz equations \cite{Lima}. The Bethe Ansatz
equations are a system of coupled and non-linear equations introduced
in the field of statistical mechanics by Bethe in 1931 \cite{Bethe}.
In \cite{Lima}, from an appropriated change of variables, the authors
managed to reduce the Bethe Ansatz equations associated to the $XXZ$
six-vertex model to a coupled system of polynomial equations and,
for the so called two-magnon sector, they succeeded in decoupling
this system of equations, so that the solutions could be written in
terms of the roots of the following self-inversive polynomials 
\begin{equation}
P_{a}(z)=\left(\omega_{a}+1\right)z^{n}-2\omega_{a}\Delta z^{n-1}-2\Delta z+\left(\omega_{a}+1\right),\label{Bethe}
\end{equation}
where $\omega_{a}=e^{2\pi ia/n}$, $1\le a\leq n$, is one of the
$n$th roots of unity and $\Delta$ is a parameter specific to the
model. In this work, Vieira and Lima-Santos studied the distribution
of the roots of the polynomials (\ref{Bethe}), and was stated that
all the roots of $P_{a}(z)$ lie on the complex unit circle $U$ and
are simple if
\begin{equation}
\left|\Delta\right|<\left|\frac{\omega_{a}+1}{2}\right|.\label{Delta1}
\end{equation}
On the other hand, if the inequality 
\begin{equation}
\left|\Delta\right|>\frac{n}{n-2}\left|\frac{\omega_{a}+1}{2}\right|\label{Delta2}
\end{equation}
is satisfied, then was stated that all but two simple roots $s$ and
$\omega_{a}/s$ of $P_{a}(z)$ are on $U$ (we remark that these results
follow directly from the theorems presented above). Furthermore, it
was verified as well that $P_{a}(z)$ may present multiple roots on
$U$ whenever any of the inequalities (\ref{Delta1}) or (\ref{Delta2})
are replaced by an equality. Finally, we might see that the polynomial
$P_{a}(z)$ becomes a Salem polynomial whenever $\Delta$ is a half-integer
greater than $1$ $-$ the appearance of Salem polynomials in the
solutions of the Bethe Ansatz equations was quite surprising, since
they were found only in a few fields of mathematical physics so far
\cite{Smyth,Hironaka}.

\thanks{The author thanks A. Lima-Santos for the motivation, comments and
discussions and also the anonymous referee of this paper for his valuable
suggestions. This work was supported by the São Paulo Research Foundation
(FAPESP), grant \#2012/02144-7 and grant \#2011/18729-1.}

\end{document}